\documentclass[reqno,12pt]{amsart}
\usepackage{amscd,amssymb}

\numberwithin{equation}{section}

\theoremstyle{plain}
\newtheorem{Thm}[subsection]{Theorem}
\newtheorem{Cor}[subsection]{Corollary}
\newtheorem{Lem}[subsection]{Lemma}
\newtheorem{Prop}[subsection]{Proposition}

\theoremstyle{definition}
\newtheorem{Def}[subsection]{Definition}

\theoremstyle{remark}

\newtheorem{Rem}[subsection]{Remark}

\renewcommand{\rm}{\normalshape}

\newif\ifShowLabels
\ShowLabelstrue
\newdimen\theight
\def\TeXref#1{%
        \leavevmode\vadjust{\setbox0=\hbox{{\tt
                \quad\quad  {\small \rm #1}}}%
        \theight=\ht0
        \advance\theight by \lineskip
        \kern -\theight \vbox to
        \theight{\rightline{\rlap{\box0}}%
        \vss}%
        }}%

\ShowLabelsfalse

\renewcommand{\sec}[2]{\section{#2}\label{S:#1}%
        \ifShowLabels \TeXref{{S:#1}} \fi}
\newcommand{\ssec}[2]{\subsection{#2}\label{SS:#1}%
        \ifShowLabels \TeXref{{SS:#1}} \fi}

\newcommand{\refs}[1]{Section ~\ref{S:#1}}
\newcommand{\refss}[1]{Section ~\ref{SS:#1}}
\newcommand{\reft}[1]{Theorem ~\ref{T:#1}}
\newcommand{\refl}[1]{Lemma ~\ref{L:#1}}

\newcommand{\refc}[1]{Corollary ~\ref{C:#1}}
\newcommand{\refd}[1]{Definition ~\ref{D:#1}}
\newcommand{\refr}[1]{Remark ~\ref{R:#1}}
\newcommand{\refe}[1]{\eqref{E:#1}}

\newenvironment{thm}[1]%
        { \begin{Thm} \label{T:#1}  \ifShowLabels \TeXref{T:#1} \fi }%
        { \end{Thm} }

\renewcommand{\th}[1]{\begin{thm}{#1} \sl }
\renewcommand{\eth}{\end{thm} }

\newenvironment{lemma}[1]%
        { \begin{Lem} \label{L:#1}  \ifShowLabels \TeXref{L:#1} \fi }%
        { \end{Lem} }
\newcommand{\lem}[1]{\begin{lemma}{#1} \sl}
\newcommand{\elem}{\end{lemma}}

\newenvironment{propos}[1]%
        { \begin{Prop} \label{P:#1}  \ifShowLabels \TeXref{P:#1} \fi }%
        { \end{Prop} }
\newcommand{\prop}[1]{\begin{propos}{#1}\sl }
\newcommand{\eprop}{\end{propos}}

\newenvironment{corol}[1]%
        { \begin{Cor} \label{C:#1}  \ifShowLabels \TeXref{C:#1} \fi }%
        { \end{Cor} }
\newcommand{\cor}[1]{\begin{corol}{#1} \sl }
\newcommand{\ecor}{\end{corol}}

\newenvironment{defeni}[1]%
        { \begin{Def} \label{D:#1}  \ifShowLabels \TeXref{D:#1} \fi }%
        { \end{Def} }
\newcommand{\defe}[1]{\begin{defeni}{#1} \sl }
\newcommand{\edefe}{\end{defeni}}

\newenvironment{remark}[1]%
        { \begin{Rem} \label{R:#1}  \ifShowLabels \TeXref{R:#1} \fi }%
        { \end{Rem} }
\newcommand{\rem}[1]{\begin{remark}{#1}}
\newcommand{\erem}{\end{remark}}

\newcommand{\eq}[1]%
        { \ifShowLabels \TeXref{E:#1} \fi
           \begin{equation} \label{E:#1} }
\newcommand{\eeq}{\end{equation}}

\newcommand{\prf}{ \begin{proof} }
\newcommand{\eprf}{ \end{proof} }

\newcommand\alp{\alpha}         
\newcommand\bet{\beta}

\newcommand\tet{\theta}         
\newcommand\iot{\iota}

\newcommand\ome{\omega}         

\newcommand\calL{{\mathcal{L}}}

\newcommand\calO{{\mathcal{O}}}

\newcommand\RR{\mathbb{R}}

\newcommand\ZZ{\mathbb{Z}}

\newcommand\CC{\mathbb{C}}

\newcommand\nek{,\ldots,}
\newcommand\sdp{\times \hskip -0.3em {\raise 0.3ex
\hbox{$\scriptscriptstyle |$}}} 

\newcommand\IM{\operatorname{Im}}

\newcommand\ind{\operatorname{ind}}

\newcommand\spin{\operatorname{spin}}

\newcommand\oz{{\overline{z}}}

\newcommand\tilE{{\widetilde{E}}}

\newcommand\tilp{{\widetilde{p}}}

\newcommand\tilmu{{\widetilde{\mu}}}

\newcommand{\w}{\omega}

\renewcommand{\O}{\calO}

\newcommand{\hm}[2]{{H^{#2}(M_{#1},\calO(E_{#1}))}}

\newcommand{\m}{\backslash}

\newcommand\ch{\operatorname{char}}
\newcommand\mult{\operatorname{mult}}

\newcommand{\hf}{\frac12}

\newcommand{\ka}{K\"ahler }
\begin{document}

\title{Symplectic cutting of K\"ahler  manifolds}
\author{Maxim Braverman}
\address{Institute of Mathematics\\
         The Hebrew University   \\
         Jerusalem 91904 \\
         Israel
         }
\email{maxim@math.huji.edu}
\thanks{This research was partially supported by grant No. 96-00210/1 from
the United States-Israel Binational Science Foundation (BSF), Jerusalem, Israel}
\begin{abstract}
We obtain estimates on the character of the cohomology of an $S^1$-equivariant
holomorphic vector bundle over a \ka manifold $M$ in terms of the cohomology of the
Lerman symplectic cuts and the symplectic reduction of $M$. In particular, we prove and
extend inequalities conjectured by Wu and Zhang \cite{WuZhang}.

The proof is based on constructing a flat family of complex spaces $M_t \ (t\in\CC)$
such that $M_t$ is isomorphic to $M$ for $t\not=0$, while $M_0$ is a singular reducible
complex space, whose irreducible components are the Lerman symplectic cuts.
\end{abstract}

\maketitle

\sec{introd}{Introduction}

Let $M$ be a smooth \ka manifold of complex dimension $n$ endowed with a holomorphic
Hamiltonian action of the circle group $S^1$. Let $\mu:M\to \RR$ denote the moment map
for this action. Assume that 0 is a regular value of $\mu$ and that $S^1$ acts freely
on $\mu^{-1}(0)$. Using a construction of E.~Lerman, \cite{Lerman-cut}, one can "cut"
$M$ into two smooth \ka manifolds $M_+$ and $M_-$ endowed with a holomorphic circle
action. The symplectic reduction $M_{red}=\mu^{-1}(0)/S^1$ of $M$ is embedded into
$M_\pm$ as a connected component of the fixed point set.

Let $E$ be an equivariant holomorphic vector bundle over $M$. Then $E$ induces a
holomorphic vector bundle $E_{red}$ over $M_{red}$ and equivariant holomorphic bundles
$E_\pm$ over $M_\pm$. In this paper we show that there are Morse-type inequalities
which estimate the character of the $S^1$-action on the cohomology $\hm{}{*}$  of the
sheaf of holomorphic sections of $E$ in terms of the cohomology \/ $\hm{\pm}{*}, \
\hm{red}{*}$ \/ of the sheaves of holomorphic sections of the bundles $E_\pm$ and $E_{red}$
respectively. These inequalities were conjectured by Wu and Zhang \cite{WuZhang} for
the case when $E$ is a pre-quantum line bundle.

As a consequence, we obtain a new geometric proof of the "gluing formula" for the index
of $E$ (\cite{DGMW,Meinr-GS}, see \refe{glue})   for the case when the manifold $M$ is
K\"ahler.

Our proof is based on the following geometric construction which, we believe, is
interesting by itself. We consider the union $M_{cut}$ of $M_\pm$ along $M_{red}$. Thus
$M_{cut}$ is a singular reducible complex space, whose smooth irreducible components
$M_\pm$ intersect by the symplectic reduction $M_{red}$. We show that $M_{cut}$ may be
considered as a deformation of $M$. More precisely, we construct a family  $M_t$ of
complex spaces parameterized by a complex parameter $t$, such that $M_t$ is complex
isomorphic to $M$ for any $t\not=0$ while $M_0$ is complex isomorphic to $M_{cut}$. It
turns out that $M_t$ is {\em a flat family of complex spaces}. That implies that the
dimension (and also the character) of the cohomology of $M_t$ is an upper
semi-continuous function of $t$. In particular, the character of $\hm{cut}{*}$ is
greater than the character of $\hm{}{*}$ (a partial order on the ring of characters is
introduced in \refd{polyn}). Moreover, there are Morse-type inequalities (cf.
\reft{Mcut>M}) which estimate the character of $\hm{}{*}$ in terms of the character of
$\hm{cut}{*}$.

The cohomology $\hm{cut}{*}$  of the space $M_{cut}$ can be, in turn, calculated by
means of a Mayer-Vietoris-type long exact sequence via the cohomology of $M_\pm$ and
$M_{red}$ (cf. \refss{Mcut}). That leads to estimates for the cohomology of $M$ in
terms of the cohomology of $M_\pm$ and $M_{red}$.

The paper is organized as follows. In \refs{main}, we formulate our main results. In
\refs{family}, we present our geometric construction of the family of complex spaces
and prove some important properties of this family. Finally, in \refs{proof}, we
present the proof of \reft{Mcut>M}.

\subsection*{Acknowledgments}
I would like to thank I.~Zakharevich for very useful and inspiring discussions. It was
I.~Zakharevich who suggested to consider the union of the Lerman symplectic cuts
$M_\pm$ as a singular complex space.

I would like to thank the University of Warwick, where this work was completed, for
hospitality.

\sec{main}{Main results}

In this section we formulate the main results of the paper. All these results are
consequences of \reft{Mcut>M}, which will be proved in \refs{proof}.

\ssec{char}{Weights and formal characters} Irreducible representation
of the circle group $S^1=\{e^{i\tet}:\, \tet\in\RR\}$ are classified by integer {\em
weights} (here we use the identification of the Lie algebra of $S^1$ with $\RR$ which
takes the {\em negative} primitive lattice element, $-2\pi i\in i\RR= Lie(S^1)$, to
$1$). A representation of weight $k\in \ZZ$ is isomorphic to the complex line $\CC$ on
which the element $e^{i\tet}\in S^1$ acts by multiplications by $e^{-ik\tet}$.

If $W$ is a finite dimensional representation of \/ $S^1$ \/ we denote by $\mult_k(W)$
the multiplicity of the weight $k\in\ZZ$ in $W$.

The {\em formal character} of $W$ is the formal sum
$$
        \ch(W) \ = \ \sum_{k\in\ZZ}\mult_k(W)e^{-ik\tet}.
$$
It lies in the ring \/ $\calL=\ZZ[e^{i\tet},e^{-i\tet}]$ \/ of Laurent polynomials in
$e^{i\tet}$ with integer coefficients. This ring is called the {\em ring of formal
characters} of the circle group.

\ssec{mom-red}{Momentum map and symplectic reduction}
Let $V$  denote the vector field on $M$ that generates the $S^1$-action and let $\w$
denote the \ka form on $M$. We will assume that $S^1$-action is {\em Hamiltonian}, i.e.
there is a moment map $\mu: M\to\RR$ such that $\iot_V\ome=d\mu$. Note (\cite{Frankel})
that it is always the case if the fixed-point set of $S^1$ on $M$ is non-empty.

Assume that $0\in \RR$ is a regular value of the moment map $\mu$. Then
$\mu^{-1}(0)\subset M$ is a smooth submanifold endowed with a locally free action of
$S^1$. We will assume that this action is free. Then the quotient space
$M_{red}=\mu^{-1}(0)/S^1$ is a smooth symplectic manifold called the {\em symplectic
reduction of $M$ at level $0$}. The symplectic form $\w_{red}$ on $M_{red}$ is defined
by the condition that its lift on $\mu^{-1}(0)$ coincides with the restriction of $\w$
on $\mu^{-1}(0)$.

Recall now that our manifold $M$ is \ka and that the \ka structure on $M$ is preserved
by the circle action. In this case, {\em the $S^1$ action can be canonically extended
to a holomorphic action of the group of nonzero complex numbers $\CC^*$} (cf.
\cite[Lemma~3.3]{GuiSter82}). The set
$$
        M_s \ = \ \big\{z\cdot x: \ z\in \CC^*, x\in \mu^{-1}(0)\subset M \big\},
$$
called {\em the set of stable points} for the $\CC^*$ action, is an open submanifold of
$M$, \cite[Lemma~4.5]{GuiSter82}, and the $\CC^*$ action on $M_s$ is free. Obviously,
the quotient of $M_s$ by this action is diffeomorphic to the reduced space:
\eq{M/C}
        M_{red} \ \cong \ M_s/\CC^*.
\end{equation}
The equation \refe{M/C},  defines {\em a  canonical complex structure on $M_{red}$}.
This structure is, in fact, K\"ahler, and the corresponding \ka form coincides with the
form $\w_{red}$ defined above.

Let now $E$ be a holomorphic vector bundle over $M$ which is equivariant for the $S^1$
action. Then the $\CC^*$ action on $M$ can be also lifted on $E$. There is a unique
holomorphic vector bundle $E_{red}$ over $M_{red}$ such that its pullback under the
projection \/ $M_s\to M_{red}$ \/ is isomorphic to the restriction of $E$ on $M_s$.

\ssec{cut}{Symplecting cuttings}
We now recall the Lerman symplectic cutting construction, \cite{Lerman-cut}. Let
$\CC_{\pm}$ denote the complex one-dimensional representations of the circle group of
weights $\pm 1$ respectively. We endow both $\CC_+$ and $\CC_-$ with the standard \ka
form $\ome=\frac{i}{2}d z\wedge d\oz$. The diagonal actions of $S^1$ on
$M\times\CC_\pm$ are Hamiltonian and the corresponding moment maps are
$\mu\mp\hf|z|^2$. One checks easily that 0 is a regular value for each one of these
moment maps. Let us denote by $M_\pm$  the symplectic quotients of $M\times\CC_\pm$ at
level $0$. The action of $S^1$ on the first factor of $M\times \CC_\pm$ reduces to a
Hamiltonian action on $M_\pm$. Thus, $(M_\pm,\ome_\pm)$ are smooth symplectic manifolds
with Hamiltonian $S^1$-actions. The reduced space $M_{red}$ is embedded  into $M_\pm$
as one of the connected components (still denoted by $M_{red}$) of the fixed points
set; the compliments $M_\pm\m M_{red}$ are $S^1$-equivariantly symplectomorphic to
$\mu^{-1}(\RR^\pm)\subset M$, respectively. We refer to $M_\pm$ as {\em symplectic cuts
of $M$}.

The pull-back of the bundle $E$ under the natural projection $M\times \CC_\pm\to M$ is
an equivariant vector bundle over $M\times\CC_\pm$.  Hence (cf. \refss{mom-red}), it
induces holomorphic vector bundles $E_\pm$ over $M_\pm$. One of the most important
facts about the cohomology of the symplectic cuts is the {\em gluing formula} (cf.
\cite{DGMW,Meinr-GS})
\begin{multline}\label{E:glue}
        \sum_{p=0}^n (-1)^p\ch\hm{}{p} \ = \ \sum_{p=0}^n (-1)^p\ch\hm{+}{p} \\
        \ + \ \sum_{p=0}^n (-1)^p\ch\hm{-}{p}
        \ - \ \sum_{p=0}^{n-1} (-1)^p\dim_{\CC} \hm{red}{p}.
\end{multline}

\rem{symplectic}
Though the individual cohomology $\hm{}{p}$ has sense only for complex manifold $M$,
the alternating sums which appear in \refe{glue} may be defined in the case when $M$ is
an almost complex manifold. The formula \refe{glue} remains true  for this more general
case \cite{DGMW,Meinr-GS} (see also \cite{SiKaTo} were the gluing formula is obtained
in a still more general situation).
\erem

Let us return to the situation when $M$ is K\"ahler. The aim of this paper is to
strengthen the gluing formula \refe{glue} in order to obtain an information about
individual cohomology $\hm{}{p}$ of $M$ in terms of the cohomology of $M_\pm$ and
$M_{red}$.

\ssec{sing}{Symplectig cutting as a singular space}
Both manifolds $M_\pm$ contain the symplectic reduction $M_{red}$ as a submanifold.
Consider the union
$$
        M_{cut}= M_+\cup_{M_{red}} M_-
$$
along $M_{red}$. Then $M_{cut}$ is a singular reducible complex space whose irreducible
components are $M_\pm$ and whose only singularities are the "double points" in
$M_{red}$. Let $E_{cut}$ denote the vector bundle over $M_{cut}$ whose restriction onto
$M_\pm$ is equal to $E_\pm$. The advantage of considering the singular space $M_{cut}$
rather then two disconnected manifolds $M_\pm$ is that $M_{cut}$ may be considered as a
deformation of $M$ (cf. \refs{family}). This implies (cf. \reft{Mcut>M}) estimates on
\/ $\ch H^*(M,\O(E))$  \/ in terms of the character \/ $\ch H^*(M_{cut},\O(E_{cut}))$ \/
of the cohomology of the sheaf of holomorphic sections of $E_{cut}$. The cohomology \/
$H^*(M_{cut},\O(E_{cut}))$ \/  my be, in turn, calculated in terms of the cohomology of
the sheaves \/ $E_\pm$ \/ and \/ $E_{red}$ \/ (cf. \refss{Mcut} bellow). That gives an
estimate on
\/ $\ch H^*(M,\O(E))$ \/ in terms of the spaces \/ $M_{\pm},M_{red}$.

To formulate the result we need the following
\defe{polyn} Let  $q(\tet)= \sum_{k\in\ZZ}q_ke^{-ik\tet}\in \calL$ \/  be a formal
  character of $S^1$, we say $q(\tet)\ge 0$ if $q_k\ge 0$ for all
  $k\in\ZZ$. For two characters $p,q\in \calL$, we say that $p\ge q$
  if $p-q\ge 0$.

  Let $Q(\tet,t)= \sum_{m=0}^n q_m(\tet)t^m\in \ \calL[t]$ be a
  polynomial of degree $n$ with coefficients in $\calL$, we say
  $Q(\tet,t)\ge0$ if $q_m(\tet)\ge0$ for all $m$.
\edefe

Our first result is the following Morse-type inequalities between the cohomology of $M$
and $M_{cut}$.
\th{Mcut>M} There exists a polynomial $Q(\tet,t)\in \calL[t]$, such that $Q\ge 0$ and
  \eq{Mcut>M}
        \sum_{p=0}^n t^p\ch\hm{cut}{p} \ = \ \sum_{p=0}^n t^p\ch\hm{}{p} \ + \
                        (1+t)Q(\tet,t).
  \end{equation}
\eth
\reft{Mcut>M} is proven in \refs{proof}.
\rem{Mcut>M} \ 1. \ The Morse-type inequalities \refe{Mcut>M}  imply
  $$
        \ch\hm{cut}{p} \ \ge \ \ch\hm{}{p} \quad \mbox{for any}\quad p=0\nek n.
  $$

  2. \ Substituting $t=-1$ into \refe{Mcut>M} we obtain the following  index formula
  \eq{Mcut=M}
        \sum_{p=0}^n (-1)^p\ch\hm{}{p} \ = \ \sum_{p=0}^n (-1)^p\ch\hm{cut}{p}.
  \end{equation}
\erem

\ssec{Mcut}{Cohomology of $M_{cut}$} To calculate the cohomology of $M_{cut}$
with coefficients in $\O(E_{cut})$ consider the equivariant short exact sequence of
sheaves
$$\begin{CD}
        0\to \O(E_{cut}) \ @>\alp>> \ \O(E_+)\oplus \O(E_-) \
                @>\bet>> \O(M_{red}) \ \to 0,
\end{CD}
$$
Here the map \/ $\alp$ \/ sends the section $s$ of \/ $\O(E_{cut})$ \/ to the pair \/
$(s|_{M_+},s|_{M_-})$ \/ and the map \/ $\bet$ \/ sends the pair of sections \/
$(s_+,s_-)\in \O(E_+)\oplus \O(E_-)$ \/ to the section \/
$s_+|_{M_{red}}-s_-|_{M_{red}}\in \O(M_{red})$.

By standard arguments, the above short sequence leads to an equivariant long exact
sequence in cohomology
\begin{multline}\label{E:HMcut}
        \cdots\to H^p(M_{cut},\O(E_{cut}))\to H^p(M_+,\O(E_+))\oplus H^p(M_-,\O(E_-)) \\
                \to H^p(M_{red},\O(E_{red}))\to H^{p+1}(M_{cut},\O(E_{cut}))\to\cdots
\end{multline}
We think about $M_{cut}$ as being glued from $M_\pm$ along $M_{red}$. So we refer to
\refe{HMcut} as  Mayer-Vietoris-type sequence.
\rem{WuZhang}
  Wu and Zhang \cite[Remark~4.10]{WuZhang} conjectured that, if $E$ is a {\em pre-quantum
  line bundle}, then  (for a proper choice of the moment map) the cohomology $\hm{}{p}$ may
  be calculated by a Mayer-Vietoris-type exact sequence of type \refe{HMcut}. In our terms,
  that would mean that, in this case, $\hm{}{p}$ is isomorphic to $\hm{cut}{p}$.
\erem

The long exact sequence \refe{HMcut} leads to  the following Morse-type inequalities
\begin{multline}\label{E:MorseMcut}
        \sum_{p=0}^n t^p\ch\hm{+}{p} \\
                + \sum_{p=0}^n t^p\ch\hm{-}{p} \ + \
                        \sum_{p=0}^{n-1} t^{p+1}\dim\hm{red}{p} \\
        = \ \sum_{p=0}^n t^p\ch\hm{cut}{p} \ + \ (1+t)Q(\tet,t)
\end{multline}
for some $Q(\tet,t)\ge0$. Combining with \reft{Mcut>M} we obtain the following estimate
on $\hm{}{*}$ in terms of the cohomology of $M_\pm$ and $M_{red}$
\th{morse} There exists a polynomial $Q'(\tet,t)\in \calL[t]$, such that $Q\ge 0$ and
  \begin{multline}\label{E:MorseM}
        \sum_{p=0}^n t^p\ch\hm{+}{p} \ +  \ \sum_{p=0}^n t^p\ch\hm{-}{p} \\
        \ + \
                \sum_{p=0}^{n-1} t^{p+1}\dim\hm{red}{p}
        = \ \sum_{p=0}^n t^p\ch\hm{}{p} \ + \ (1+t)Q'(\tet,t)
  \end{multline}
\eth
In the case when $E$ is a pre-quantum line bundle \reft{morse} was conjectured by Wu
and Zhang \cite[Remark~4.10]{WuZhang}.

\rem{gluing} \
1. \ The inequalities \refe{MorseMcut} a far from being exact. Hence, \reft{Mcut>M}
together with the Mayer-Vietoris sequence \refe{HMcut} give much more information about
the cohomology $\hm{}{*}$ than \reft{morse}.

2. \ The simplest consequence of \reft{morse} is the inequalities
\eq{simple}
 \begin{aligned}
  \ch\hm{+}{0} \ &+ \ \ch\hm{-}{0} \ \ge \ \ch\hm{}{0}; \\
  \ch\hm{+}{p} \ &+ \ \ch\hm{-}{p} \ + \ \dim\hm{red}{p-1} \\
   &\ge \ \ch\hm{}{p}, \quad\ \mbox{for any} \ \quad p=1\nek n.
 \end{aligned}
\end{equation}

3. \ Substituting $t=-1$ into \refe{MorseM} and using \refe{Mcut=M}, we get the gluing
formula \refe{glue}. So we obtain a new proof of the gluing formula for
\ka manifolds, which is based on the geometric construction described in \refs{family}.
Note that the standard proof of the gluing formula, \cite{DGMW,Meinr-GS,SiKaTo}, uses
the Atiyah-Segal-Singer equivariant index theorem. Dietmar Salamon pointed out that the
gluing formula for general symplectic manifold can also be proved using our geometric
construction by a method similar to \cite[Appendix~A]{McDSal2}.
\erem
\rem{combine}
In the situation considered in this paper many other Morse-type inequalities may be
obtained (cf. \cite{WuZhang,TianZhang1,Br-HM}). It would be very interesting to compare
those inequalities.
\erem

\ssec{example}{Example}
We finish this section with a very simple but typical example illustrating
\reft{morse}.

Let $M=\CC P^1$. We identify $M$ with the 2-dimensional sphere $S^2\subset\RR^3$ and we
let $S^1$ act on $M$ by rotations around the $z$-axis. This action has two fixed points
$P$ and $Q$ (the poles of the sphere). We normalize the \ka structure on $M$ and the
moment map $\mu$ so that $\mu(P)=1, \ \mu(Q)=-1$. Then the image of $\mu$ is the
interval $[-1,1]$ and all the internal points of this interval are regular values of
$\mu$.

Let $E$ be an equivariant line bundle over $M$. Then $S^1$ acts on the fibers of this
bundle over the fixed points $P$ and $Q$. Denote by $r_Q, r_P$ the weights of these
actions. It is well known that $E$ is defined up to an equivariant isomorphism by these
weights.  In particular (cf., for example, \cite[p.~330]{Witten84}) the character of
the representation of $S^1$ on the cohomology $\hm{}{p}$ is given by
\footnote{Note that our signs in the definition of weights  and characters (cf.
\refss{char}) are different from \cite{Witten84} but agree with
\cite{WuZhang,SiKaTo,Br-HM}.}
\eq{cohom}
  \begin{aligned}
        \ch \hm{}{0} &=
             \begin{cases}
                \sum_{m=r_Q}^{r_P}e^{-im\tet}, \ \ \quad&\mbox{if}\quad r_Q\le r_P;\\
                0, \ \ \quad&\text{if}\quad r_Q> r_P;
             \end{cases}  \\
        \ch \hm{}{1} &=
             \begin{cases}
                0, \quad&\text{if}\quad r_Q\le r_P;\\
                \sum_{m=r_P-1}^{r_Q-1}e^{-im\tet}, \quad&\mbox{if}\quad r_Q> r_P.
             \end{cases}
   \end{aligned}
\end{equation}
These formulas allow us to calculate the right hand side of \refe{MorseM}. Let us
calculate the left hand side of \refe{MorseM}.

Since, $M_{red}$ is a point, \/ $\dim\hm{red}{0}=1$ \/ and \/ $\dim\hm{red}{1}=0$.

Both manifolds $M_+$ and $M_-$ are isomorphic to $\CC P^1$. The weight of the fiber of
the bundle $E_+$ over $P$ is still equal to $r_P$ while the weight of the fiber over
$M_{red}$ is equal to zero. Similarly, the weight of the fiber of the bundle $E_-$ over
$Q$ is equal to $r_Q$ while the weight of the fiber over $M_{red}$ is equal to zero.

Using \refe{cohom}, one can now verify \reft{morse} in this simple case. For example,
if $r_Q=r_P=r>0$ (this corresponds to the trivial bundle $E=M\times\CC$ with a
nontrivial circle action), then the left hand side of \refe{MorseM} is equal to
$$
        \sum_{m=0}^re^{-im\tet} + t\sum_{m=0}^{r-1}e^{-im\tet}
$$
while
$$
        \sum_{p=0}^n t^p \ch\hm{}{p} = e^{-ir\tet}.
$$
It follows that in this case the polynomial $Q'$ of \reft{morse} does not depend on
$t$ and equals to \/ $\sum_{m=0}^{r-1}e^{-im\tet}$.

Note also, that {\em the inequalities \refe{MorseM}, \refe{simple} become equalities if
and only if $r_Q\le0$ and $r_P\ge0$}. This verifies the conjecture of Wu and Zhang
\cite{WuZhang} (cf. \refr{WuZhang}) in the case $M=\CC P^1$.


\sec{family}{The geometric construction}

In this section we present a geometric construction of  a complex manifolds $\Phi$  and
a holomorphic map $p:\Phi\to \CC$ such that $p^{-1}(t)$ is isomorphic to $M$ for
$t\not=0$ while $p^{-1}(0)=M_{cut}$.

\ssec{family}{}
The idea of the construction is the same as in \refss{cut}. In fact, we just combine
the constructions of $M_+$ and $M_-$ together and  consider the diagonal action of
$S^1$ on $M\times\CC_+\times\CC_-$.  The moment map for this action is given by
$$
        \tilmu(x,z_+,z_-) \ = \ \mu(x) \ - \ \hf|z_+|^2+\hf|z_-|^2,   \qquad
                        x\in M, \ z_\pm\in\CC_{\pm}.
$$
Zero is a regular value of $\tilmu$ and we define $\Phi= \tilmu^{-1}(0)/S^1$ to be the
symplectic reduction of $M\times\CC_+\times\CC_-$ at zero level.

Clearly, the  map  $\tilp: \, M\times\CC_+\times\CC_-\to \CC$  defined by the formula
$$
        \tilp: \ (x,z_+,z_-) \ \mapsto \ z_+z_-, \qquad   x\in M, \ z_\pm\in\CC_{\pm}.
$$
is $S^1$ invariant and, hence, descends to a map $p:\Phi\to \CC$. We think about $\Phi$
as family of complex manifolds $M_t=p^{-1}(t)$ parameterized by a complex parameter
$t\in\CC$.

We endow $\Phi$ with an $S^1$-action induced by the action of the circle on the first
factor of $M\times\CC_+\times\CC_-$. Note that this action preserves the fibers
$M_t=p^{-1}(t)$ of $p$. Thus $M_t \ (t\in\CC)$ are also endowed with a holomorphic
circle action.

As in \refss{cut} the bundle $E$ induces a bundle $\tilE$ over $\Phi$. We denote the
restriction of $\tilE$ on $M_t$ by $E_t$.

The following lemmas describe the fibers of the projection $p$.
\lem{MtS}
  For any $t\not=0$, the fiber $M_t=p^{-1}(t)$ is a smooth  manifold which is
  equivariantly symplectomorphic to $M$.
\elem
\prf
  Fix a nonzero number $t\in\CC$. For any $x\in M$ set
  $$
        r(x) \ = \ \sqrt{\mu(x)+\sqrt{\mu(x)^2+|t|^2}}
  $$
  and define an embedding
  \eq{it}
        i_t: \, M\to \tilmu^{-1}(0)\cap\tilp^{-1}(t) \ \subset M\times\CC_+\times\CC_-,
        \qquad
        i_t: \ x \ \mapsto \Big(x,r(x),\frac{t}{r(x)}\Big).
  \end{equation}
  Clearly, the composition $q\circ i_t$ of the above embedding with the natural projection
  $q: \, \tilmu^{-1}(0)\to \Phi$ is an equivariant diffeomorphism $M\to M_t=p^{-1}(t)$. Since
  the map $i_t:M\to M\times\CC_+\times\CC_-$ is symplectic, so is the composition $q\circ i_t$.
\eprf
\lem{MtC}
  For any $t\not=0$, the fiber $M_t=p^{-1}(t)$ is a smooth \ka manifold which is
  equivariantly complex isomorphic to $M$. In other words, there exists an
  equivariant biholomorphic map $\phi_t:M\overset{~}{\to} M_t$.

  The pullback
  $\phi^*_tE_t$ of the bundle $E_t=\tilE|_{M_t}$ is equivariantly isomorphic to $E$.
\elem
\rem{isom}
  The manifolds $M$ and $M_t$ are not isomorphic as \ka manifolds. In particular, the map
  $\phi_t$ of \refl{MtC} is different from the symplectomorphism of \refl{MtS}.
\erem
\prf
  Recall from \refss{mom-red} that the action of the circle group on $M$ extends
  canonically to a holomorphic action of  $\CC^*$ and consider the diagonal action of
  $\CC^*$ on $M\times\CC_+\times\CC_-$ (here $z\in \CC^*$ acts on the second factor by
  multiplication by $z$ and on the third factor by multiplication by $1/z$).

  Let $U\subset M\times\CC_+\times\CC_-$ denote the set of stable points for this action
  (cf. \refss{mom-red}). Then $\tilp^{-1}(t)\subset U$ for any $t\not=0$.
  Indeed, if $v\in \tilp^{-1}(t)$ and if the absolute value of the number $z\in\CC^*$
  is large enough, then $\tilmu(z\cdot v) < 0$, while $\tilmu(\frac1z\cdot v) > 0$.
  Hence, one can find $z'\in\CC^*$ such that $z'\in \tilmu^{-1}(0)$.

  Fix $t\not=0$ and consider the complex map
  $$
        j_t:M\to M\times\CC_+\times\CC_-, \qquad j_t: x \ \mapsto \big(x,1,t\big).
  $$
  Clearly, the image of $j_t$ belongs to $\tilp^{-1}(t)$ and, by the previous paragraph,
  it belongs also to the set $U$ of stable  points for $\CC^*$ action. Hence, the
  composition of $j_t$ with the quotient map $q:U\to U/\CC^*$ defines an equivariant
  holomorphic map $\phi_t:M\to M_t$. Clearly, $\phi_t$ is injective. We claim that
  $\phi_t$ is an isomorphism. By \refl{MtS}, it
  suffice to show that the image of $j_t$ contains the image of the map \refe{it}. But
  this follows from the obvious inclusion \/ $\frac1r\cdot(x,r,\frac{t}r)\in \IM(j_t)$.
  The first statement of the lemma is proven.

  Consider the commutative diagram
  \eq{diagr}
    \begin{CD}
        M\times\CC_+\times\CC_- &&\ &\ \hookleftarrow&\quad &U\cap \tilp^{-1}(t)\\
        @V{\pi}VV               \ &{\nearrow}&\quad         &@VV{q}V\\
        M                       &\ & @>\phi_t>>\quad                     &M_t
    \end{CD}
  \end{equation}
  By definition of the bundle $E_t$ we have $q^*E_t=\pi^*E|_{U\cap \tilp^{-1}(t)}$. Hence,
  using $\phi_t=q\circ j_t$ and  $\pi\circ j_t=id$, we obtain
  $$
        \phi_t^*E_t \ = \ j_t^*q^*E_t \ = \ \phi_t^*\pi^*E \ = \ E.
  $$
\eprf

\lem{M0}
  The fiber $p^{-1}(0)$ is equivariantly complex isomorphic to the space $M_{cut}$. If we
  identify $p^{-1}(0)$ with $M_{cut}$ using this isomorphism then the restriction $E_0$
  of $\tilE$ to $p^{-1}(0)$ is isomorphic to $E_{cut}$.
\elem
\prf
  The lemma is an obvious consequence of the equality
  $$
        \tilp^{-1}(0) \ = \ (M\times\CC_+\times\{0\}) \, \cup \,
                        (M\times\{0\}\times\CC_-).
  $$
\eprf

For us the most important is the following consequence of the above lemmas
\cor{Phi}
  The cohomology $H^*(M_t,\O(E_t))$ of the sheaf of holomorphic
  sections of the bundle $E_t$ is equivariantly isomorphic to $\hm{}{*}$ if $t\not=0$ and
  is equivariantly isomorphic to $\hm{cut}{*}$ if $t=0$.
\ecor


\sec{proof}{Flat morphisms. Proof of \reft{Mcut>M}}

We are in a position now to prove \reft{Mcut>M}. The proof is based on the properties
of flat morphisms in complex analysis.

\ssec{flat}{Flat morphisms}
First, we recall some basic facts about flat
morphisms. For the details we refer the reader to
\cite[Sections~II.2,III.4]{GraPetRem}.

If  \/ $X$ \/  is a complex space and  \/ $x\in X$ \/  we denote by  \/ $\O(X)$ \/  the
sheaf of holomorphic functions on  \/ $X$ \/  and by  \/ $\O_x(X)$ \/  the ring of
germs of holomorphic functions at  \/ $x$. Let  \/ $f:X\to Y$ \/  be a holomorphic map
of complex spaces. For any  \/ $y\in Y$, we denote by  \/ $X_y=f^{-1}(y)$ \/  the fiber
of
\/ $f$ \/  over  \/ $y$.

A holomorphic map \/ $f:X\to Y$ \/ is called {\em flat at a point \/ $x\in X$} \/ if
\/ $\O_x(X)$ \/ is a flat module over  \/ $\O_{f(x)}(Y)$. Here  \/ $\O_x(X)$ \/  is considered as
 \/ $\O_{f(x)}(Y)$-module via the canonical map  \/ $f^*:\O_{f(x)}(Y)\to
\O_x(X)$. The map  \/ $f$ \/  is called {\em flat} if it is flat at any point  \/ $x\in X$.

Let  \/ $f:X\to Y$ \/  be a morphism of complex spaces and suppose  \/ $V$ \/  is a
holomorphic vector bundle over  \/ $X$. For any \/ $y\in Y$, let  \/ $V_y=V|_{X_y}$ \/
denote the restriction of \/ $V$ \/  on the fiber  \/ $X_y$ \/  and let  \/ $\O(V_y)$
\/ denote the locally free sheaf of holomorphic sections of  \/ $V_y$. If  \/ $f$ \/
is a flat morphism, then, for any  \/ $y\in Y$ \/ and for any pint \/ $\eta\in Y$ \/
closed to $y$,  there exists a polynomial  \/ $Q(t)$ \/  with nonnegative integer
coefficients, such that
\eq{family}
        \sum_p t^p\dim H^p(X_y,\O(V_y)) \ = \
                \sum_p t^p\dim H^p(X_\eta,\O(V_\eta))
        \ + \ (1+t)Q(t).
\end{equation}

The equation \refe{family} implies, in particular,  that the function  \/ $y\mapsto
\dim H^p(X_y,\O(V_y))$ \/  is upper semi-continuous, while the {\em holomorphic index} of
the fibers
$$
        \ind(V_y) \ = \ \sum_p (-1)^p\dim H^p(X_y,\O(V_y))
$$
is locally constant on  \/ $Y$.

Suppose that in the situation described above a compact Lie group  \/ $G$ \/  acts
holomorphically on  \/ $X$ \/  and  \/ $Y$ \/  and that this action commutes with  \/
$f$. If the vector bundle  \/ $E$ \/  is equivariant with respect to this action then
\/ $G$ \/  acts on the cohomology of the fibers. Let  \/ $\ch H^p(X_y,\O(V_y)), \ y\in Y$ \/
denote the character of this action. Then, for any point  \/ $\eta\in Y$ \/  closed
enough to  \/ $y$, there exists a polynomial  \/ $Q(t,\tet)\in \calL[t]$ \/  (cf.
\refd{polyn}) such that
 \/ $Q\ge 0$ \/  and
\eq{eq-family}
        \sum_p t^p\ch H^p(X_y,\O(V_y)) \ = \
                \sum_p t^p\ch H^p(X_\eta,\O(V_\eta))
        \ + \ (1+t)Q(t,\tet).
\end{equation}

\ssec{proof}{Proof of \reft{Mcut>M}}
It follows now from \refe{eq-family} and \refc{Phi}, that in order to prove
\reft{Mcut>M} it suffices to show that the projection  \/ $p:\Phi\to \CC$ is flat. But,
by a theorem of Kaup and Kerner, \cite[Ch.~II, Theorem~2.13]{GraPetRem}, {\em any open
holomorphic map of smooth complex manifolds is flat}. Since,  \/ $p$ \/  is open the
theorem is proven.
\hfill $\square$

\providecommand{\bysame}{\leavevmode\hbox to3em{\hrulefill}\thinspace}

\end{document}